
\catcode`\@=11  


%
%
\hsize = 15.8truecm
\vsize = 24.2truecm

\parskip=\bigskipamount

\input amssym.def
\input amssym.tex

\font\tenbmit   = cmmib10
\font\sevenbmit = cmmib7
\font\fivebmit  = cmmib5

\newfam\bmitfam

\textfont\bmitfam     = \tenbmit
\scriptfont\bmitfam   = \sevenbmit
\scriptscriptfont\bmitfam = \fivebmit

\def\leq{\leqslant}	\def\geq{\geqslant}

\def\frac#1#2{{\textstyle {#1\over #2}}}	\def\half{\frac{1}{2}}
\def\quarter{\frac{1}{4}}

\def\C{{\Bbb C}}	\def\D{{\Bbb D}}	
	\def\P{{\Bbb P}}	
\def\R{{\Bbb R}}

\def\A{{\Bbb A}}

\newsymbol\square 1003
\def\endproof{\hfill$\square$\bigskip}

\def\epsilon{\varepsilon}



\outer\long\def\proclaim #1 #2 #3
{
    \medbreak\noindent{\bf #1} \quad #2\par
    \noindent{\sl #3}\par\penalty-100
}

\outer\def\proof{\noindent{\sl Proof:}\par\nobreak}

\outer\def\beginsection#1\par
{   
    \vskip 0pt plus 0.5\vsize \penalty-250
    \vskip 0pt plus -0.5\vsize \bigskip\bigskip\vskip\parskip
    \message{#1}\leftline{\bf #1}\nobreak\smallskip
}

\outer\def\beginsubsection#1\par
{
    \vskip 0pt plus 0.1\vsize \penalty-250
    \vskip 0pt plus -0.1\vsize \bigskip\bigskip\vskip\parskip
    \message{#1}\leftline{\bf #1}\nobreak\smallskip
}


\input epsf.tex

\catcode`\@=12



\centerline{\bf SPHERICAL AND HYPERBOLIC LENGTHS OF IMAGES OF ARCS}
\bigskip
\centerline{T. K. CARNE}
\bigskip

\centerline{Abstract}

Let $f : \D \to \C$ be an analytic function on the unit disc which is in the
Dirichlet class, so the Euclidean area of the image, counting multiplicity, is
finite.  The Euclidean length of a radial arc of hyperbolic length $\rho$ is
then $o(\rho^{1/2})$.  In this note we consider the corresponding results when
$f$ maps into the unit disc with the hyperbolic metric or the Riemann sphere
with the spherical metric.  Similar but not identical results hold.

\bigskip

\centerline{1. {\it Introduction}}

Let $f : \D \to \C$ be an analytic function on the unit disc $\D$ which is in
the Dirichlet class, so the Euclidean area of the image $f(\D)$, counting
multiplicity, is finite.  Keogh [K] showed that a radial arc $[0, z]$ in the
disc is mapped to an arc with Euclidean length ${\cal E}$ that satisfies
$$
{\cal E} = o(\rho(0, z)^{1/2})
$$
where $\rho$ is the hyperbolic length in $\D$.  The paper [BC] explored the
result in some detail.  In this paper we will consider the corresponding
results when the image domain is either the disc or the extended complex
plane.  For these we will use the natural Riemannian metrics: the hyperbolic
metric on the unit disc and the spherical metric on the extended complex
plane.  Similar results hold in these cases, essentially because we can
localise the result to a small disc with an Euclidean image.  However, there
are interesting differences and the arguments make the importance of the
hyperbolic metric still more apparent.

We will consider domains $A$ with a Riemannian metric $ds = \lambda_A(z)
|dz|$.  Here $\lambda_A$ is a strictly positive function on $A$ giving the
density of the metric.  If $f : A \to B$ is an analytic map, then the
derivative has norm
$$
||f'(z_o)||_{A \to B} = |f'(z_o)| {{\lambda_B(f(z_o))}\over{\lambda_A(z_o)}}
\ .
$$
This is the factor by which $f$ changes infinitesimal lengths at the point
$z_o$ for the metrics on $A$ and $B$.  The area is changed by the square of
this factor.

On the complex plane $\C$ we will use the Euclidean metric $|dz|$ with density
$1$.  On the unit disc $\D$ we will use the hyperbolic metric with density 
$$
\lambda_H(z) = {2\over{1-|z|^2}}\ .
$$
This has constant curvature $-1$.  On the extended complex plane we will use
the spherical metric with density
$$
\lambda_S(z) = {2\over{1+|z|^2}}\ .
$$
This has constant curvature $+1$ and is isometric, under stereographic
projection, with the unit sphere in Euclidean $\R^3$.  This is the Riemann
sphere and we will denote it by $\P$.  The subscripts $H, E, S$ will be used
to specify the hyperbolic, Euclidean and spherical metrics respectively.  So
$d_H, d_E, d_S$ are the distances on these three spaces and $\A_H, \A_E, \A_S$
are the area measures.

Let $f : \D \to B$ be an analytic map into some domain $B$ with density
$\lambda$.  The area of the image, counting multiplicity is 
$$
\int_0^1 \int_0^{2\pi} |f'(re^{i\theta})|^2 \lambda(f(re^{i\theta}))^2 \;
d\theta\; r\, dr  =
\int_\D ||f'(z)||_{E\to E}^2 \; d\A_E(z) \ .
$$
We will abbreviate this to $\A_B(f(\D))$ when there is no chance of confusion.
It is more natural to write this as an area integral over the unit disc using
the hyperbolic metric:
$$
\A_B(f(\D)) = \int_\D ||f'(z)||_{H \to B}^2 \; d\A_H(z)\ .
$$

Let $f : \D \to \C$ be an analytic map with $\A_E(f(\D)) < \infty$, so $f$ is
in the Dirichlet class.  In [BC] it was shown that 
$$
||f'(z)||_{H\to E} \leq \left({{\A_E(f(\D))}\over{4\pi}}\right)^{1/2}\ .
\eqno{(1.1)}
$$
Choose a direction from the origin and consider a radial arc of hyperbolic
length $\rho$ from $0$ in this direction.  The image of this arc will have
Euclidean length $L_E(\rho)$.  By integrating $(1.1)$ along the radial arc we
can show that
$$
L_E(\rho) = o(\rho^{1/2})
\eqno{(1.2)}
$$
as $\rho \nearrow \infty$, that is as the radial arc extends to the boundary.  

The purpose of this note is to consider the corresponding results for maps 
$f : \D \to \D$ or $f : \D \to \P$.  For the hyperbolic image, both $(1.1)$
and $(1.2)$ hold.  In the spherical case, $(1.1)$ fails when the area
$\A_S(f(\D))$ is large.  However, it does hold when this area is small and
this is sufficient to establish $(1.2)$.  

In all cases, the arguments are very similar in spirit to those in [BC] and
so they are not laboured.  That paper considers many analogues and extensions
of the results and these, similarly, can be established for hyperbolic and
spherical images.

\bigskip

\centerline{2. {\it The norm of the derivative}}

We wish to establish bounds on the derivative $||f'(z)||_{H\to B}$ in terms of
the area $\A_B(f(\D))$.  In the Euclidean case this is very simple and is
already done in Keogh's paper [K].  We include a proof to compare with later
results for the hyperbolic and spherical images.

\proclaim{Proposition 2.1} {}
{Let $f : \D \to \C$ be an analytic map with the Euclidean area of the image
$\A_E(f(\D))$, counting multiplicity, finite.  Then 
$$
||f'(z_o)||_{H\to E} \leq \left({{\A_E(f(\D))}\over{4\pi}}\right)^{1/2}\ .
$$
}
\proof
Composing $f$ with a hyperbolic isometry will not alter the area of the
image.  Hence we can assume that $z_o = 0$.  Let $f(z) = \sum a_n z^n$ be the
power series for $f$.  Then $||f'(0)||_{H\to E} = |f'(0)|/2 = |a_1|/2$.
The usual integration of the corresponding Fourier
series gives 
$$
\A_E(f(\D)) = \pi \sum_{n=1}^\infty n |a_n|^2 \ .
$$
So we certainly have 
$$
||f'(0)||_{H\to E}^2 = {{|a_1|^2}\over{4}} \leq {{\A_E(f(\D))}\over{4\pi}}
$$
as required.  The inclusion map $f : \D \hookrightarrow \C$ shows that the
constant $1/4\pi$ in this result is the best possible.
\endproof

We can also apply this result to analytic maps $f : \Delta \to \C$ defined on
any hyperbolic domain $\Delta$.  In particular, let $\Delta$ be the ball of
hyperbolic radius $\rho_o$ about $z_o$:
$$
\Delta = B(z_o, \rho_o) = \{ z \in \D : d_H(z_o, z) < \rho_o \}\ .
$$
Let $r_o = \tanh \half \rho_o$ (so the disc of hyperbolic radius $\rho_o$
centred on $0$ has Euclidean radius $r_o$).  Then the M\"obius transformation 
$$
T : z \mapsto {{r_o z+z_o}\over{1 +\overline{z_o} r_o z}}
$$
maps the unit disc conformally onto $\Delta$.  So we can apply the Proposition
to $g = f \circ T$ to obtain
$$
r_o ||f'(z_o)||_{H\to E} \leq
\left({{\A_E(f(\Delta))}\over{4\pi}}\right)^{1/2}
\ .
$$
This gives us a way of localising the result and hence of applying it when the
image domain has a different metric.
 
\proclaim{Proposition 2.2} {}
{Let $f : \D \to \D$ be an analytic map with the hyperbolic area of the image
$\A_H(f(\D))$, counting multiplicity, finite.  Then 
$$
||f'(z_o)||_{H\to H} \leq \left({{\A_H(f(\D))}\over{4\pi}}\right)^{1/2}\ .
$$
}
\proof
By composing with hyperbolic isometries, we may assume that $z_o = 0$ and
$f(z_o) = 0$.  Then $||f'(z_o)||_{H \to H} = |f'(0)|$.

The hyperbolic density in $\D$ is always at least $2$, so $\A_E(f(\Delta))
\leq \frac{1}{4} \A_H(f(\Delta))$.  Therefore, Proposition 2.1 shows that 
$$
||f'(0)||_{H\to H} = 2 ||f'(0)||_{H \to E} \leq
2 \left({{\A_E(f(\D))}\over{4\pi}}\right)^{1/2} \leq
\left({{\A_H(f(\D))}\over{4\pi}}\right)^{1/2}
$$
as required.
\endproof

For the map $f : z \mapsto \epsilon z$ we have $||f'(0)||_{H \to H} =
\epsilon$ while $\A_H(f(\D)) = 4\pi \epsilon^2/(1-\epsilon^2)$.  So we see
that the constant $1/4\pi$ in the proposition is the best possible.

The corresponding result for the spherical metric fails.  For consider a
univalent map $k : \D \to \P$ with $k(0) = 0$ and whose image has spherical
area $4\pi = \A_S(\P)$.  For example, the Koebe function $k : z\mapsto
z/(1-z)^2$.  This has $\A_S(k(\D)) = 4\pi$ and
$$
||k'(0)||_{H\to S} = |k'(0)| \not= 0\ .
$$
For any $\lambda \not= 0$, the map $f(z) = \lambda k(z)$ also has $\A_S(k(\D))
= 4\pi$ but $||f'(0)||_{H\to S} = |\lambda|$.  So we can not have any
inequality of the form
$$
||f'(0)||_{H\to S}^2 \leq c \A_S(f(\D))\ .
$$
Nonetheless, the result does hold provided that the spherical area of the
image is sufficiently small.

\proclaim{Proposition 2.3} {}
{Let $f : \D \to \P$ be an analytic map with the spherical area of the image
$\A_S(f(\D))$, counting multiplicity, less than $2\pi$.  Then 
$$
||f'(z_o)||_{H\to S} \leq c \A_E(f(\D))^{1/2}
$$
for some constant $c$.
}
\proof
We may assume that $z_o = 0$ and $f(z_o) = 0$.  So $||f'(0)||_{H\to S} =
|f'(0)|$.  

Choose $\delta > 0$ so that the spherical area of the ball  
$$
B_S(c, \delta) = \{ z\in \P : d_S(c, z) < \delta \}
$$
is less than $\half \pi$.  For this $\delta$, find a maximal set of points
$c_0, c_1, c_2, \dots , c_K$ with $c_0 =0$, the other points $c_1, c_2, c_K$
in $\P\setminus f(\D)$, and all of the distances $d_S(c_i, c_j) \geq \delta$
for $i\not= j$.  Since $K$ is maximal, there can be no point of $\P\setminus
f(\D)$ lying outside the balls $B_S(c_j, \delta)$.  Hence these $K+1$ balls
cover all of $\P\setminus f(\D)$ and so
$$
(K+1) \A_S(B_S(0, \delta)) \geq \A_S(\P\setminus f(\D)) \geq 2\pi\ .
$$
This implies that 
$$
K+1 \geq {{2\pi}\over{\A_S(B(0, \delta))}} > 4\ .
$$
So we can find $3$ points $w_0, w_1, w_\infty$ in $\P\setminus f(\D)$ which
satisfy
$$
d_S(0, w_i) > \delta\ ,\ \ d_S(w_i, w_j) > \delta
\eqno{(2.1)}
$$
for $i\not= j$.  We can now compare the spherical metric on $\P$ with the
hyperbolic metric on $\P\setminus \{w_0, w_1, w_\infty\}$.

The three punctured sphere $\P\setminus \{w_0, w_1, w_\infty\}$ has a
hyperbolic metric and the conditions $(2.1)$ allow us to estimate its
properties uniformly.  From now on, consider $\P\setminus \{w_0, w_1,
w_\infty\}$ with this metric.  This hyperbolic metric can be written as
$\lambda(z) ds_S(z)$ where $ds_S$ denotes the infinitesimal spherical metric
and $\lambda$ is the density relative to this spherical metric.  The function
$\lambda$ is bounded away from $0$ and tends to $+\infty$ at the three
punctures.  The conditions $(2.1)$ imply that there are constant $K, K', K'' >
0$, depending only on $\delta$, with the following properties:

\itemitem{(a)} $\lambda(0) \geq K$;

\itemitem{(b)} the spherical balls $B_S(w_i, K')$ are
disjoint and at least hyperbolic distance $1$ from $0$ in $\P\setminus \{w_0,
w_1, w_\infty\}$; 

\itemitem{(c)} The hyperbolic density $\lambda(z) \leq K''$ for all points $z$
within a hyperbolic distance $1$ from $0$ in $\P\setminus \{w_0,
w_1, w_\infty\}$.  
 
(Compare this with the estimates in [A], 1-9. There is a M\"obius
transformation $T : \P \to \P$ that maps our three punctured sphere
$\P\setminus \{w_0, w_1, w_\infty\}$ onto the standard three punctured sphere
$\P\setminus \{0, 1, \infty\}$.  The conditions $(2.1)$ show that this
transformation only distorts the metrics by a controlled amount.)

Now $f : \D \to \P\setminus \{w_0, w_1, w_\infty \}$ is an analytic map
between two hyperbolic domains so the Schwarz -- Pick lemma implies that $f$
is a contraction for the hyperbolic metrics.  This implies that $f$ maps the
hyperbolic disc $\Delta = B_H(0, 1)$ with hyperbolic radius $1$ into the
region 
$$
\{ w\in \P\setminus \{w_0, w_1, w_\infty \} : d_H(0, w) < 1 \}\ .
$$
Condition $(a)$ above shows that 
$||f'(0)||_{H\to S} \leq {2\over K} ||f'(0)||_{H\to H}$.
Condition $(c)$ shows that the hyperbolic area $\A_H(f(\Delta))$ is at most
$K''^2$ times the spherical area $\A_S(f(\Delta))$.  Finally, we can apply
Proposition 2.2 to the map $f|_\Delta : \Delta \to \P\setminus \{w_0, w_1,
w_\infty \}$ (or, more properly, to its lift to the universal cover of
$\P\setminus \{w_0, w_1, w_\infty \}$).  This gives 
$$
r_o ||f'(0)||_{H\to H} \leq \left({{\A_H(f(\Delta))}\over{4\pi}}\right)^{1/2}
\leq
\left( {{K''^2 \A_S(f(\Delta))}\over{4\pi}} \right)^{1/2}
$$
for $r_o = \tanh \half$.  Putting all of these together gives 
$$
||f'(0)||_{H\to S} \leq {{2K''}\over{r_o K}} 
\left({{\A_S(f(\D))}\over{4\pi}}\right)^{1/2}
$$
as required. 
\endproof

The argument used in the proof shows that, for any $C < \A_S(\P) =
4\pi$, there is a constant $K(C)$ with
$$
||f'(z_o)||_{H\to S} \leq K(C) \A_S(f(\D))^{1/2}
$$
provided that $\A_S(f(\D)) \leq C$.

\bigskip

\centerline{3. {\it The lengths of image arcs}}

We will need to introduce some notation that will apply throughout the
remainder of the paper.

Let $f : \D \to B$ be an analytic map into one of the domains $B = \D, \C,
\P$.  Let $\gamma : [0, \infty) \to \D$ be a radial hyperbolic geodesic with
unit speed and argument $\theta$, so $\gamma(t) = r e^{i\theta}$ where 
$r = \tanh \half t$.  Consider the arc $\gamma[0, \rho_o]$ with hyperbolic
length $\rho_o$.  This has an image of length
$$
L_B(\rho_o) = \int_0^{\rho_o} ||f'(\gamma(t)||_{H\to B} \; dt
$$
in the metric for $B$.  

For each $\rho < \infty$, the image of the hyperbolic ball $B_H(0, \rho) =
\{z\in \D : d_H(0, z) < \rho \}$ has finite area, which we denote by $A(\rho)
= \A_B(f(\Delta))$.   Then $A(\rho)$ is an increasing function.  If the area
of the entire image $\A_B(f(\D))$ is finite, then $A(\rho)$ is bounded above
and converges to $\A_B(f(\D))$ as $\rho \nearrow \infty$.

We wish to apply the propositions of \S 2 to the function $f$ restricted to
hyperbolic discs $B_H(\gamma(t), \delta)$ for some hyperbolic radius
$\delta$. This disc lies inside the ball $B_H(0, t+\delta)$ and outside the
disc $B_H(0, t-\delta)$.  So the area $\A_B(f(\Delta)) \leq A(t+\delta) -
A(t-\delta)$. 

Initially, let us consider the Euclidean case $f : \D \to \C$.  For this
Proposition 2.1 gives us
$$
||f'(\gamma(t))||_{H\to E} \leq {1\over{\tanh \half \delta}} 
\left({{A(t+\delta) - A(t-\delta)}\over{4\pi}}\right)^{1/2}\ .
$$
So integrating and applying the Cauchy -- Schwarz inequality gives
$$
\eqalign{
L_E(\rho_o) &= \int_0^{\rho_0} ||f'(\gamma(t))||_{H\to E} \; dt \cr
&\leq
\int_0^{\rho_0} {1\over{\tanh \half \delta}} 
\left({{A(t+\delta) - A(t-\delta)}\over{4\pi}}\right)^{1/2} \; dt \cr
&\leq
\left( \int_0^{\rho_0} {1\over{4\pi\tanh^2 \half \delta}} 
\left(A(t+\delta) - A(t-\delta)\right) \; dt\right)^{1/2} 
\left( \int_0^{\rho_0} 1 \; dt\right)^{1/2} \cr
&\leq
\left( \int_0^{\rho_0} {1\over{4\pi \tanh^2 \half \delta}} 
\left(A(t+\delta) - A(t-\delta)\right) \; dt\right)^{1/2} 
\rho_o^{1/2} \cr
}
$$
Suppose that the area $\A_E(f(\D))$ is finite.  Then $A(t)$ increases to its
limiting value $\A_E(f(\D))$.  So the integral 
$$
\int_0^{\rho_0} {1\over{\tanh^2 \half \delta}} 
\left({{A(t+\delta) - A(t-\delta)}\over{4\pi}}\right) \; dt
$$
is bounded independently of $\rho_o$.  Hence $L_E(\rho_o) = O(\rho_o^{1/2})$.
More carefully, we can apply the above result to the arc from some value $u_o$
up to $\rho_o$.  The integrand 
$
A(t+\delta) - A(t-\delta)
$
is then no more than 
$
\A_E(f(\D)) - A(u_o - \delta)
$.  This tends to $0$ as $u_o \nearrow \infty$, so we see that $L_E(\rho) =
o(\rho_o^{1/2})$ as $\rho_0 \to \infty$.  Thus we have reproved Keogh's
Theorem as in [BC].

\proclaim{Theorem 3.1} {(Keogh)}
{Let $f : \D \to \C$ be an analytic map with $\A_E(f(\D))$ finite.  For any
argument $\theta$, the Euclidean length $L_E(\rho_o)$ of the image under $f$
of the radial arc $[0, r_o e^{i\theta}]$ of hyperbolic length $\rho_0$ (so
$r_o = \tanh \half \rho_o$) satisfies 
$$
L_E(\rho_o) = o(\rho_o^{1/2}) \qquad \hbox{ as } \rho_o \to \infty.
$$
}
\endproof

Note that the proof given above only required the application of Proposition
2.1 to functions where the area of the image is small, for we knew that $A(t)
\nearrow \A_E(f(\D))$ as $t\nearrow \infty$.  Hence Propositions 2.2 and 2.3
give the corresponding results for maps into the hyperbolic plane and the
Riemann sphere.

\proclaim{Theorem 3.2} {}
{Let $f : \D \to \D$ be an analytic map with $\A_H(f(\D))$ finite.  For any
argument $\theta$, the hyperbolic length $L_H(\rho_o)$ of the image under $f$
of the radial arc $[0, r_o e^{i\theta}]$ of hyperbolic length $\rho_0$
satisfies
$$
L_H(\rho_o) = o(\rho_o^{1/2}) \qquad \hbox{ as } \rho_o \to \infty.
$$
}

\proclaim{Theorem 3.3} {}
{Let $f : \D \to \P$ be an analytic map with $\A_S(f(\D))$ finite.  For any
argument $\theta$, the spherical length $L_S(\rho_o)$ of the image under $f$
of the radial arc $[0, r_o e^{i\theta}]$ of hyperbolic length $\rho_0$
satisfies
$$
L_S(\rho_o) = o(\rho_o^{1/2}) \qquad \hbox{ as } \rho_o \to \infty.
$$
}

In [BC] examples were constructed of functions $f : \D \to \C$ showing that
the power $\half$ in Theorem 2.1 is the best possible.  Since these examples
had $f$ bounded and the hyperbolic, Euclidean and spherical metrics are
Lipschitz equivalent on any compact region inside unit disc, these examples
also show that the power is best possible in Theorems 3.2 and 3.3.

It is easy to adapt the argument used above to more general situations.
Suppose, for example, that the area $A(t)$ grows to infinity but we have
control on the rate of growth.  Then we can use the same ideas to obtain a
bound $L(\rho_o) = O(\rho_o^\alpha)$ for suitable exponents $\alpha$.  

In more detail, fix $\alpha > 1$.  A different splitting of the integrands in
the appeal to the Cauchy -- Schwarz inequality above gives us
$$
\eqalign{
L_E(\rho_o) 
&\leq
\int_0^{\rho_0} {1\over{\tanh \half \delta}} 
\left({{A(t+\delta) - A(t-\delta)}\over{4\pi}}\right)^{1/2} \; dt \cr
&\leq
\left( \int_0^{\rho_0} {1\over{4\pi t^{\alpha - 1} \tanh^2 \half \delta}} 
\left(A(t+\delta) - A(t-\delta)\right) \; dt\right)^{1/2} 
\left( \int_0^{\rho_0} t^{\alpha - 1} \; dt\right)^{1/2} \cr
&\leq
\left( \int_0^{\rho_0} {1\over{4\pi t^{\alpha - 1}\tanh^2 \half \delta}} 
\left(A(t+\delta) - A(t-\delta)\right) \; dt\right)^{1/2} 
\left({{\rho_o^{\alpha}}\over{\alpha}}\right)^{1/2} \cr
}
$$
For the first integrand, we can ignore the behaviour near $0$ and write
$$
\eqalign{
\int^{\rho_0} &{1\over{t^{\alpha - 1}\tanh^2 \half \delta}} 
\left(A(t+\delta) - A(t-\delta)\right) \; dt
= \cr
&\qquad = \int^{\rho_0-\delta} {1\over{(t-\delta)^{\alpha - 1}\tanh^2 \half \delta}} 
A(t) \; dt 
- 
\int^{\rho_0+\delta} {1\over{(t+\delta)^{\alpha - 1}\tanh^2 \half \delta}} 
A(t) \; dt \ .\cr}
$$
So we want the integral 
$$
\int^\infty 
\left({1\over{(t-\delta)^{\alpha - 1}}} - {{1\over{(t+\delta)^{\alpha - 1}}}}
\right) 
{{A(t)}\over{\tanh^2 \half \delta}} \; dt
$$
to converge.  Then it will follows that $L_E(\rho_o) = o(\rho_o^{\alpha/2})$
as $\rho_0 \nearrow \infty$.

The mean value theorem immediately gives 
$$
{1\over{(t-\delta)^{\alpha - 1}}} - {{1\over{(t+\delta)^{\alpha - 1}}}} 
\leq
2 \delta (\alpha - 1) (t - \delta)^{-\alpha}\ .
$$
So we see that we require the area $A(t)$ to grow sufficiently slowly that 
$$
\int^\infty {\delta\over{\tanh \half\delta}} {{A(t)}\over{(t-\delta)^\alpha}}
\; dt
$$
converges.

Using Proposition 2.2 in place of Proposition 2.1 gives us the same results for
functions into the unit disc.  In order to use Proposition 2.3 to obtain
corresponding results for meromorphic functions into the Riemann sphere we
need to ensure that $A(t+\delta) - A(t-\delta)$ is sufficiently small for 2.3
to apply. In order to achieve this the radius $\delta$ usually needs to
decrease as $t$ increases.

The natural class of functions to consider here is those of finite order, so
the Nevanlinna characteristic $T(r)$ satisfies 
$$
\limsup_{r \to 1} {{T(r)}\over{\log {1\over{1-r}}}} \ < \infty\ .
$$
(See [T].)
This implies that 
$$
\int^\infty T(r) (1-r)^{k-1} \; dr < \infty
$$
for some $k$.  The derivative $T'(r)$ is $S(r)/r$ where 
$$
S(r) = {{\A_S(f(\{z : |z| < r \}))}\over{4\pi}} = {{A_S(t)}\over{4\pi}}
$$
for $t = \log (1+r)/(1-r)$, which is the hyperbolic radius corresponding to
the Euclidean radius $r$.  So, integrating by parts gives 
$$
\int^\infty A_S(t) \exp (-(k+1)t) \; dt \ < \ \infty\ .
$$
In order to obtain estimates for these functions we would need to take the
radius $\delta$ tending exponentially to $0$ as $t$ increased to $\infty$.
The details do not seem inspiring.

\bigskip

\centerline{4. {\it Examples}}

As in [BC], it is useful to consider examples that limit what can happen to
the lengths of the images of arcs.  As there we can construct many examples
defined on a region
$$
S = \{ x+iy \in \C : |y| < h(x) \}
$$
for some slowly increasing function $h$.  This will be a hyperbolic
simply-connected domain, so it is conformally equivalent to the unit disc.
Let $q : \D \to S$ be the conformal map which fixes the origin.
The hyperbolic metric is Lipschitz equivalent to the pseudo-hyperbolic metric
which has density $1/d_E(z, \partial S)$.  The positive real axis is then a
hyperbolic geodesic in $S$ and the length of the segment from the origin to
$t_o$ is approximately
$$
\rho_o = \int_0^{t_o} {1\over{h(t)}} \; dt\ .
$$
The hyperbolic disc $B(0, \rho_o)$ is then certainly contained in the part of
$S$ to the left of $\{x+iy \in S : y< t_o\}$.  

In [BC] the function $q$ itself was considered.  For us, we need to follow $q$
by a mapping that is an isometry from the positive real axis with the
Euclidean metric to the hyperbolic plane or the Riemann sphere with their
metrics.  For this we follow $q$ by an exponential map.  This gives us
corresponding examples for maps into the hyperbolic plane or the Riemann
sphere.

It may also be worth considering analogous examples for maps into the Riemann
sphere where we constrain the Nevanlinna characteristic rather than the
spherical area of the image.  If the analytic map $f : \D \to \P$ has
$\A(f(\D))$ finite, then we certainly have 
$$
T(r) = T(\half) + \int_\half^r {{\A(f(\D))}\over{4\pi s}} \; ds =  
T(\half) + {{\A(f(\D))}\over{4\pi}} \log 2r\ .
$$
So the Nevanlinna characteristic is bounded.  However, we do not have
$L_S(\rho_o) = o(\rho_o^{1/2})$ for every function $f$ with bounded
characteristic.  

Consider the universal cover of an annulus $\{ z\in \C : R^{-1} < |z| < R \}$,
say $q : \R^2_+ \to \{ z\in \C : R^{-1} < |z| < R \}$ defined on the upper
half-plane $\R^2_+$.  We can arrange for the hyperbolic geodesic $\{iy : y> 0
\}$ to be mapped to the unit circle.  If the segment from $i$ to $Ki$ is
mapped to one complete circuit of the circle, then so are the segments from
$K^n i$ to $K^{n+1} i$ for every integer $n$.  Consequently we see that the
hyperbolic geodesic from $i$ to $i\infty$ has an image with $L_S(\rho_o) \sim
\rho_o$ and not a half power.

An almost identical argument gives the same conclusion for the Blaschke
product 
$$
B(z) = \prod_{n = -\infty}^{-1} (-1)\left({{2^n i - z}\over{2^n i + z}}\right)
\prod_{n=0}^\infty \left({{2^n i - z}\over{2^n i + z}}\right)
$$
with zeros evenly spaced hyperbolically along the imaginary axis.  The image
of the positive imaginary axis is now the curve that traces out repeatedly the
line segment between the two critical values of $B$.  

This shows that we can not hope for a better inequality than $L_S(\rho_o) =
O(\rho_o)$ for functions with bounded Nevanlinna characteristic.  Even this is
untrue, as the following example shows.

Let $(y_n)$ be a strictly increasing sequence of strictly positive real
numbers with $\sum 1/y_n$ convergent.  Then there is a Blaschke product 
$$
B(z) = \prod \left({{iy_n - z}\over{iy_n + z}}\right)
$$
with zeros at the points $i y_n$.  This product is symmetric about the
imaginary axis with 
$$
B(-\overline{z}) = \overline{B(z)}\ .
$$
So, between any two successive zeros $i y_n$ and $i y_{n+1}$, there is a
single critical point.  We will be interested in the case where the $y_n$
converge slowly to $\infty$, so the hyperbolic distances $d_H(i y_n, i
y_{n+1}) = \log y_{n+1}/y_n$ decrease to $0$ as $n\to \infty$.  The Blaschke
product is a contraction for the hyperbolic metric because of the Schwarz --
Pick lemma, so the images under $f$ of the critical points must converge to
$0$.  The image of the positive imaginary axis then traces out line segments on
the real axis between successive critical values and converges to $0$.  Now
consider the geodesic $\gamma(t) = -1 + i e^t$ in the upper half-plane.  As
$t \to \infty$, this becomes closer, in the hyperbolic metric, to the positive
imaginary axis.  Hence we see that image $B(\gamma)$ traces out a path that
never takes the value $0$ and winds (negatively) about $0$.  Since
$d_H(\gamma(t), i e^t) < 1/e^t$, we see that, for large $n$, the path
$B(\gamma[iy_n, iy_{n+1}])$ completes approximately half a circuit about $0$.  

Define $f : \R^2_+ \to \P$ by
$$
f(z) =   {{B(z+1)}\over{B(z-1)}}\ .
$$
This is an analytic function with bounded Nevanlinna characteristic because of
Fatou's theorem, which says that the ratio of two bounded analytic functions
has bounded characteristic (see [N] or Proposition 4.1 later).  For a point $i
y$ on the imaginary axis we have
$$
f(iy) = {{B(1 + iy)}\over{B(-1 +iy)}} = {{B(1 + iy)}\over{\
\overline{B(1+iy)}\ }}
\ .
$$
So the path $f(\gamma[iy_n, iy_{n+1}])$ lies on the unit circle and completes
approximately one circuit about $0$ for each sufficiently large $n$.  If we
set $\rho_o = \log y_n$ to be the hyperbolic distance from $i$ to $i y_n$,
then we have 
$$
L_S(\rho_o) = L_S(\log y_n) \sim 2\pi n\ .
$$

For example, take $y_n = n^2$.  Then 
$$
L_S(\rho_o) \sim 2\pi \exp \half \rho_o\ .
$$
So we certainly do not have $L_S(\rho_o) = O(\rho_0)$.

The issue here is that the Nevanlinna characteristic depends crucially on the
choice of the origin.  When we try to apply the arguments of Theorem 3.3 we
take each point of the radial arc as the origin and so we require a bound on
the Nevanlinna characteristic independent of the position of the origin.  We
will see that such a stronger, uniform condition is enough to give
$L_P(\rho_o) = O(\rho_o)$.  It will be useful in this context to have a
precise form of Fatou's Theorem:

\proclaim{Proposition 4.1} {Fatou}
{For a meromorphic function $f : \D \to \P$ the following
conditions are equivalent.

\item{(a)} $f$ has bounded characteristic with $T(1) \leq K$;

\item{(b)} There are two analytic functions $f_0, f_\infty : \D \to
\C$ with $f = f_0/f_\infty$, 
$$
|f_0(z)|^2 + |f_\infty(z)|^2 \leq 1 \ \hbox{ and }\ 
|f_0(0)|^2 + |f_\infty(0)|^2 \geq e^{-2K} 
$$
for all $z\in \D$.
}
\proof
First we will consider this result when the function $f$ extends
analytically across the boundary $\D$ with no zeros or poles on
$\{z : |z|=1 \}$.  We can then obtain the general result by applying
this to the functions restricted to discs of radius $r < 1$.

We will need to use the chordal distance $k(w, w')$ on the Riemann sphere.
This is the length of the chord in $\R^3$ joining the points $w, w'
\in \P$, so
$$
k(w, w') = {{2|w-w'|}\over{\sqrt{1+|w|^2} \sqrt{1+|w'|^2}}}
\qquad \hbox{ and } \qquad
k(w, w') = 2 \sin \half d_S(w, w')\ .
$$

Write $f(z) = \displaystyle{{{B_0(z)}\over{B_\infty(z)}} \exp h(z)}$
where $B_0, B_\infty$ are finite Blaschke products on the zeros and
poles of $f$ and $h$ is analytic on the closed unit disc.  Let $u_0,
u_\infty$ be continuous functions on the closed disc, harmonic on the
interior and with boundary values
$$
u_0(\zeta) = \log k(f(\zeta), 0) \ ,\ \ 
u_\infty(\zeta) = \log k(f(\zeta), \infty)\ \ 
\hbox{ for } \ |\zeta| = 1 .
\eqno{(4.1)}
$$
Then
$$
u_0(\zeta) - u_\infty(\zeta) = \log {{k(f(\zeta), 0)}\over{k(f(\zeta),
\infty)}} = 
\log |f(\zeta)| = {\rm Re\,} h(\zeta)
$$
and so $u_0(z) - u_\infty(z) = {\rm Re\,} h(z)$ for each $z\in \D$.  We
can choose harmonic conjugates $\tilde{u}_o$ and $\tilde{u}_\infty$
with $(u_0 + i\tilde{u}_0) - (u_\infty + i\tilde{u}_\infty) = h$.  Now
set
$$
f_0 = \half B_0 \exp (u_0 + i\tilde{u}_0) \ \hbox{ and }\ 
f_\infty = \half B_\infty \exp (u_\infty + i\tilde{u}_\infty) .
$$
This certainly gives $f = f_0/f_\infty$.  For $|\zeta| = 1$ we have
$$
\eqalign{
|f_0(\zeta)|^2 + |f_\infty(\zeta)|^2 &= 
\quarter \exp 2u_0(\zeta) \ +\ \quarter \exp 2u_\infty(\zeta) \cr
&=
\quarter \left( k(f(\zeta), 0)^2 + k(f(\zeta), \infty)^2 \right) .\cr
}
$$
The points $0, \infty$ and $f(z)$ are the vertices of a right-angle
triangle in ${\Bbb P}$, so Pythagoras' theorem shows that 
$k(f(z), 0)^2 + k(f(z), \infty)^2 = 2^2$.  Hence
$$
|f_0(\zeta)|^2 + |f_\infty(\zeta)|^2 = 1 \qquad\hbox{ for } \ |\zeta|
 = 1 .
$$
The First Nevanlinna Theorem shows that
$$
T(1) = N(1; 0) + m(1; 0) = 
-\log |B_0(0)| + \int_0^{2\pi} \log {{k(f(0),
0)}\over{k(f(e^{i\theta}), 0)}} \; {{d\theta}\over{2\pi}} 
$$
and, since $u_0$ is harmonic, this is
$$
T(1) =  -\log |B_0(0)| - u_0(0) + \log k(f(0), 0) =
-\log 2|f_0(0)| + \log k(f(0), 0) .
$$
Therefore, $|f_0(0)| = \half k(f(0), 0) \exp -T(1)$ and,
similarly, \hfil\break 
$|f_\infty(0)| = \half k(f(0), \infty) \exp -T(1)$.
Hence,
$$
|f_0(0)|^2 + |f_\infty(0)|^2 = \quarter \left(k(f(0), 0)^2 + k(f(0),
 \infty)^2 \right) \exp -2 T(1) = \exp -2 T(1) .
$$
This shows that (a) implies (b).  It is a little simpler to reverse
this argument to prove that (b) implies (a).

Finally, for any analytic function $f : {\Bbb D} \to {\Bbb P}$ we can
apply the above result to 
$$
f_r : {\Bbb D} \to {\Bbb P} \ ;\ \ z\mapsto f(rz)
$$
for those $r < 1$ with no zeros or poles of $f$ on $\{z : |z|=r \}$.
By taking locally uniform limits as $r\nearrow 1$ we obtain the
proposition in general.
\endproof

We will say that an analytic function $f : \D \to \P$ has {\it uniformly
bounded characteristic} if there is a constant $C$ with each of the functions
$$
f_{z_o} : z\mapsto f\left({{z+z_o}\over{1+\overline{z_o}z}}\right)
$$
having characteristic $T(f_{z_o}; r)$ at most $C$ for every $r < 1$.  This is
saying that the Nevanlinna characteristic is bounded independently of the
origin $z_o$ we choose in $\D$.  The smallest value for $C$ is clearly
invariant under composing $f$ with hyperbolic isometries of $\D$.  Fatou's
theorem immediately gives us:

\proclaim{Corollary 4.2} {}
{A meromorphic function $f : {\Bbb D} \to {\Bbb P}$ has uniformly
bounded characteristic if and only if $f = f_0/f_\infty$ for two bounded analytic
functions $f_0, f_\infty$ with 
$$
\delta < (|f_0|^2 + |f_\infty|^2)^{1/2} \leq 1
$$
for some $\delta > 0$.
}
This means that the pair of functions $f_0, f_\infty$ are Corona data.

\proof
Suppose that $f : \D \to \P$ has bounded characteristic.  Then we
can write $f = f_0/f_\infty$ where
$$
f_0 = \half B_0 \, \exp (u_0 + i \tilde{u}_0)\  \ ;\ \  
f_\infty = \half B_\infty \exp (u_\infty + i\tilde{u}_\infty) 
$$
and $u_0, u_\infty$ are the harmonic functions with boundary values $(4.1)$.
Theorem 4.1 shows that
$$
\exp -2T(f_{z_o}; 1) = |f_0(z_0)|^2 + |f_\infty(z_o)|^2 .
$$
Hence, $|f_0(z_o)|^2 + |f_\infty(z_o)|^2 \geq \delta^2$ if and only if
$T(f_{z_o}; 1) \leq - \log \delta$.
\endproof

We can now prove that the inequality $L_S(\rho_o) = O(\rho_o)$ does hold for
functions $f$ with uniformly bounded Nevanlinna characteristic.

\proclaim{Theorem 4.3} {}
{Let $f : \D \to \P$ be an analytic function with uniformly bounded Nevanlinna
characteristic.  Then 
$$
L_S(\rho_o) = O(\rho_o)\ .
$$
}
\proof
The corollary shows that we can write $f$ as $f_0/f_\infty$ where $f_0,
f_\infty$ are both bounded analytic functions.  Moreover the function
$$
F : \D \to \C^2 \ ;\ \ z \mapsto (f_0(z), f_\infty(z))
$$
will satisfy $\delta \leq ||F(z)|| \leq 1$ for the Euclidean norm $||\ ||$ on
$\C^2$.  

The Schwarz -- Pick lemma, applied to $f_0$ and $f_\infty$, shows that 
$$
||F'(z)|| (1- |z|^2) \leq 2\ .
$$
Also, a simple calculation gives
$$
\eqalign{
||f'(z)||_{H\to S} &= 
\left({{2|f'(z)|}\over{1 + |f(z)|^2}}\right) \left({{1-|z|^2}\over{2}}\right)
=
{{|f_0'(z) f_\infty(z) - f_0(z) f_\infty'(z)|}\over{|f_0(z)|^2 +
|f_\infty(z)|^2}}
(1-|z|^2) \cr
&\leq 
{{||F'(z)|| \, ||F(z)||}\over{||F(z)||^2}} (1-|z|^2) 
\leq 
{{||F'(z)||}\over{||F(z)||}} (1-|z|^2) \ .\cr}
$$
So we have 
$$
||f'(z)||_{H\to S} \leq
{{2}\over{||F(z)||}} 
\leq {2\over\delta}\ .
$$
Integrating this along a hyperbolic geodesic give the result.
\endproof

\bigskip
\centerline{\it References}

\item{{\bf [A]}} 
L.V. Ahlfors, 
{\sl Conformal invariants, topics in geometric function theory},
McGraw-Hill, 1973.

\item{{\bf [BC]}} 
A.F. Beardon and T.K. Carne, 
Euclidean and hyperbolic lengths of images of arcs,
{\it Proc. London Math. Soc.}, (2007) [to appear].

\item{{\bf [K]}} 
F.R. Keogh, 
A property of bounded schlicht functions,
{\sl J. London Math. Soc.}, 29 (1954), 379-382.

\item{{\bf [N]}}
R. Nevanlinna, {\it Eindeutige analytische Funktionen}, 
Springer Verlag, Berlin , 1936. 

\item{{\bf [T]}} 
M. Tsuji,
{\sl Potential theory in modern function theory},
Maruzen, Tokyo, 1959.

\vskip 3mm
{\obeylines \parindent = 0pt
\parskip = 0pt
Department of Pure Mathematics and Mathematical Statistics,
Centre for Mathematical Sciences,
Wilberforce Road,
Cambridge. CB3 0WB
UK. 
}

\bigskip
\noindent
{\it tkc@dpmms.cam.ac.uk}

\bye